\documentstyle[12pt]{article}
\input amssym.def
\input amssym.tex

\begin{document}

\thispagestyle{empty}

\pagestyle{myheadings}

\markright{\centerline{TORRES}}

\def\bbN{{\rm I}\!{\rm N}}
\def\bbR{{\rm I}\!{\rm R}}
\def\bbZ{{\rm Z}\!\!{\rm Z}}
\def\bbP{{\rm I}\!{\rm P}}
\def\bbH{{\rm I}\!{\rm H}}
\def\bbG{{\rm G}\!{\rm G}}

\begin{center}
\null\vskip2.5truecm

{\bf ON CERTAIN N--SHEETED COVERINGS\\
OF CURVES AND NUMERICAL SEMIGROUPS WHICH\\
 CANNOT BE REALIZED AS WEIERSTRASS SEMIGROUPS\\}
\vskip1truecm
{Fernando Torres \footnote{\normalsize This work was iniciated while the
author was in Rio de Janeiro (IMPA) with a grant from Cnpq -- Brazil.}\\}
{\small International Centre for Theoretical Physics}\\
{\small Maths. Group, P.O. Box 586 - 34100 Trieste-Italy}\\
{\small feto@ictp.trieste.it}
\end{center}
\vskip0.5truecm
\begin{abstract}
A curve $X$ is said to be of type $(N,\gamma)$ if it is an $N$--sheeted
covering of a curve of genus $\gamma$ with at least one totally ramified
point. A numerical semigroup $H$ is said to be of type $(N,\gamma)$ if it
has $\gamma$  positive  multiples of $N$ in $[N,2N\gamma]$ such that its
$\gamma^{th}$ element is $2N\gamma $ and $(2\gamma+1)N \in H$.

If the genus of $X$ is large enough and $N$ is prime, $X$ is of
type $(N,\gamma)$ if and
only if there is a point $P \in X$ such that the Weierstrass semigroup at
$P$ is of type $(N,\gamma)$ (this generalizes the case of double
coverings of curves). Using the proof of this result and the Buchweitz's
semigroup, we can construct numerical semigroups
that cannot be realized as
Weierstrass semigroups
although they might satisfy  Buchweitz's criterion.
\end{abstract}
\vskip0.5truecm
\section{Introduction.}

In Weierstrass Point Theory one associates a numerical semigroup to any
non--singular point $P$ of a
projective, irreducible, algebraic curve defined over an
algebraically closed field. This semigroup is called the
Weierstrass semigroup at $P$ and is the same for all but finitely many
points. These finitely many points, where exceptional values of the
semigroup occur, are called the Weierstrass points of the curve. They
carry a lot of information about the curve.

In 1893, Hurwitz asked about the characterization of the numerical
semigroups
which arise as Weierstrass semigroups. See
[E-H] for further historical information. Long after that, in 1980
Buchweitz [B1] showed that not every numerical semigroup can occur as a
Weierstrass semigroup, but  has to satisfy
the following criterion (which can be extended to singular
curves by [S1, p. 124]),
\medskip

{\bf (BC):} ``Let $P$ be a non--singular point of a projective,
irreducible, algebraic curve defined over an algebraically closed field.
If $n \ge 2$  and $g$ is the arithmetical genus of the curve, then
the cardinality of the set of sums of $n$ gaps $G_n$ at $P$
 is bounded above by the dimension
of the pluri--canonical divisor $nC_{X}$ which is $ (2n-1)(g-1)$ ".
\medskip

Moreover, Buchweitz showed  that for every
integer
$n \ge 2$ there exist semigroups which do not satisfy the above
criterion. However, as was noticed by Oliveira [O, Thm. 1.5] and
Oliveira-St\"ohr [O-S, Thm. 1.1], this criterion cannot be applied to
symmetric and
quasi-symmetric semigroups because, in the first case $\#G_n = (2n-1)(g-1)$
and in the second  case $\#G_n = (2n-1)(g-1) - (n-2)$.
\medskip

Let $X$ be a projective, irreducible, non--singular algebraic curve of genus
$g$
defined over an algebraically closed field of characteristic $p\ge 0$.
For a point $P$ of $X$, let $H(P)$ denote the Weierstrass semigroup at $P$.
A curve is called  $\gamma$--hyperelliptic if it  is a double covering of a
curve of genus $\gamma$. A numerical semigroup $H$ will be called
$\gamma$--hyperelliptic if it has $\gamma$  positive even elements in the
interval $[2,4\gamma]$ such that its $\gamma^{th}$ element is
$4\gamma$, and, $4\gamma+2
\in H$. The reason for this terminology is the following result that has
been proved for $p =0$ in [T, Thm. A and Remark 3.10]:
\medskip

{\bf Theorem 1.} If $g\ge
6\gamma +4$, then the following statements are equivalent:
\begin{list}
\setlength{\rightmargin 0cm}{\leftmargin 0cm}
\itemsep=0.5pt
\item[(i)] $X$ is $\gamma$--hyperelliptic.

\item[(ii)] There exists $P\in X$ such that $H(P)$ is
$\gamma$--hyperelliptic.

\item[(iii)] For some integer $2\le i \le \gamma+2$ such that $2\gamma +i
\not\equiv 0$ (mod $3$) if $ i < \gamma$, there
exists a base--point--free linear
system on $X$ of projective dimension $\gamma +i$ and degree $4\gamma+2i$.
\end{list}

Using
the proof of the implication (ii)
$\Rightarrow$ (i) above and  Buchweitz's example, St\"ohr
[T,Scholium 3.5] obtained symmetric numerical semigroups which cannot be
realized as Weierstrass semigroups at points of non--singular curves. Since
 St\"ohr's examples does not depend on the characteristic $p$,
 we have that these examples cannot even be realized as
Weierstrass semigroups for non-singular curves defined in positive
characteristic. However, symmetric semigroups can be realized as
Weierstrass semigroups of Gorenstein curves [S]. Let $B_{H}$ be the
semigroup ring
over an algebraically closed field of characteristic zero associated to
a numerical semigroup $H$. By the work of Pinkham [P], $B_{H}$ is
smoothable if and only if $H$ is a Weierstrass semigroup. On the other
hand, $H$ is symmetric if and only if $B_{H}$ is Gorenstein (cf.
[Ku]). Hence St\"ohr's examples also show that the Gorenstein condition
does not imply smoothability.\medskip

The aim of the present paper is to give a generalization of the case
$p=0$ of Theorem 1 to certain covers of degree $N$. This will allow us to
construct  numerical semigroups with a given last gap (the result of
Theorem 1 just works for symmetric semigroups) that cannot be realized
as Weierstrass semigroups at points of non--singular algebraic curves
although
they might satisfy Buchweitz's criterion above  (cf. Section 4).  We
introduce the following definitions:\medskip

{\bf Definition 1.} A curve $X$ is said to be of type $(N,\gamma )$ if it
is an $N$--sheeted covering of a curve of genus $\gamma$ with at least one
totally ramified point.
\medskip

{\bf Definition 2.} Let $N>0$, $\gamma \ge 0$ be integers.
A numerical subsemigroup $H$ of $(\bbN ,+)$ is said to
be of type $(N,\gamma )$ if the following conditions are satisfied:
\begin{list}
\setlength{\rightmargin 0cm}{\leftmargin 0cm}
\itemsep=0.5pt
\item[(a)] $H$ has $\gamma$ positive multiples of $N$ in the
interval $[N,2N\gamma]$,

\item[(b)] the $\gamma^{th}$ element of $H$ is $2N\gamma$, and,

\item[(c)] $(2\gamma +1)N\in H$.
\end{list}

These  definitions are related to each other by Lemma 3.4 below. Let
$X$ be a curve of type $(N,\gamma)$ and $P$ be a totally ramified
point of $X$. If
$p=char(K) \not\!\vert\ N$ or if $N$ is prime and the genus of $X$ is large
enough, then the semigroup
$H(P)$
satisfies conditions (a), (c) above, and $2N\gamma \in H(P)$.
Moreover, $H(P)$ will be of type $(N,\gamma)$ whenever
the genus is also large enough (Corollary 3.5) (the appearance of these
semigroups  in such a
 context justifies Definition 2).
In the case where
the covered curve is given as a quotient of $X$ by an
automorphism, the above results are included in an implicit way in T. Katos's
paper
[K, p.395]. From an arithmetical point of view, we can say that every
numerical semigroup $H$ is a semigroup of type $(1,\gamma)$ with $\gamma$
equal to the genus of $H$. Moreover, given a positive integer $N$ there
exists a natural number $\gamma_N=\gamma_N(H)$ such that $H$ satisfies
conditions
(a), (c) above (with $N$ and $\gamma_N$) and $2N\gamma_N \in H$ (Lemma 2.3
(ii)). The point is that it
  does not necessarily fulfil
condition (b) (Remark 3.11 (i)).
 On the other
hand, using a result of R\"ohrl [R\"o, Th.3.1] we can construct curves that
are not
of type $(N,\gamma )$. It follows also that there exist curves $X$
of type $(N,\gamma )$ which do not cover a
curve of the form
$X/\langle T\rangle$ for $T$  an automorphism of order $N$ defined on $X$.
\medskip

We now state the main result of this paper. For $A , u, N$ and $\gamma$
integers we  define:
\begin{equation}
\rho_1(A,N,\gamma) = {A(N-1)N\over 2} + N\gamma -N +1,
\end{equation}
\begin{equation}
\rho_2(N,\gamma) = N(2N-1)\gamma - (N-1)(N+2),
\end{equation}
\begin{equation}
\rho_3(N,\gamma) = (2N-1)(N\gamma + N-1),
\end{equation}
\begin{equation}
\rho_4(A,u,N,\gamma) = {1\over 2}(N-u-1)[(A-\gamma-1)(N+u) - 2(N\gamma
+ N-1)] + \rho_3(N,\gamma).
\end{equation}

{\bf Theorem A.} Consider the following statements:
\begin{list}
\setlength{\rightmargin 0cm}{\leftmargin 0cm}
\itemsep=0.5pt
\item[(i)] $X$ is a curve of type $(N,\gamma )$.

\item[(ii)] There exists $P\in X$ such that $H(P)$ is a semigroup of type
$(N,\gamma )$.

\item[(iii)] There exists $P\in X$ and an integer $A$ such that the
linear system
$|AP|$ is base--point--free of dimension $A - \gamma$.
\end{list}
\vspace{-0.5pt}
Let $N$ be prime and $A\ge 2\gamma +1$ an integer.
\begin{list}
\setlenght{\rightmargin 0cm}{\leftmargin 0cm}
\itemsep=0.5pt
\item[(A1)] If  $N \not= p = char(K)$ and $g>
\rho_1(2\gamma,N,\gamma) =
N^2\gamma - N+1$ or if $g> \rho_2(N,\gamma)$, then (i) $\Rightarrow$ (ii).

\item[(A2)] If $g > \rho_3(N,\gamma)$, then (ii)$
\Rightarrow $ (i).

\item[(A3)] If $g>
\rho_1(A,N,\gamma)$, then (ii) $\Rightarrow$ (iii).

\item[(A4)] Let $u= u(A)$ be the biggest
integer $\le {N\gamma +N-1 \over A-\gamma -1}$. If the following conditions
\begin{equation}
A \not\equiv 0\ ({\rm mod}\ t)\ {\rm for}\ t \le {AN\over A-\gamma},\
t\not= N\ {\rm and}
\end{equation}
\begin{equation}
g > \rho_4(A,u,N,\gamma)
\end{equation}
 hold, then (iii) $\Rightarrow$ (i).
\end{list}

Notice that for $N=2$, statement (iii) of Theorem 1 is equivalent to
statement (iii) of Theorem A because the $\gamma$--hyperelliptic
involution is unique (if it exists) provided $g> 4\gamma +1$ (see e.g.
[A, Lemma 5]). Definitions (3) and (4)  are derived from
Castelnuovo's
number (see Lemma 3.7 below) applied to a linear system of type
$g^{A-\gamma}_{AN}$; definition (2) arise from Castelnuovo's bound
involving subfields of the field of rational functions of the curve
(Lemma 3.2) while definition (1) comes from arithmetical reasons (Lemma
2.1).

In Section 2 we study some arithmetical properties of the semigroups of type
$(N,\gamma )$. As in the $\gamma$--hyperelliptic case ([T, Lemma 2.6]), we find
that the multiples of $N$ contained in the semigroup are
determined by properties (a) and (c) (Lemma 2.3 (i)). We also state
sufficient conditions for a numerical semigroup to be of type $(N,\gamma)$
(Corollaries 2.5 and 2.6). If $N\ge 2$ is an integer, we have a
lower bound for the elements $h$ of the semigroup such that
 $gcd(h,N) = 1$ (Lemma 2.1). This result generalizes Lemma 2.2 in [T].
Both Lemma 2.1 and Lemma 2.3 are used to obtain linear series on curves
having a point with Weiertrass semigroup of type $(N,\gamma)$.

In Section 3 we prove Theorem A and some results concerning Weierstrass
semigroups at totally ramified points. In Remark 3.11 we discuss
the sharpness of the bounds on $g$ used in the results of the paper as well
as the necessity of hypothesis (5).

In section 4, using the proof of item (A2), we  show how to
construct
numerical semigroups of type $(N,\gamma )$ that are not realized as Weierstrass
semigroups of non--singular curves. We also construct numerical semigroups
that do not satisfy
 Buchweitz's criterion for $n =2$. These examples contain
Buchweitz's
semigroup and they may be well known, but we included them here in
order to compare them with the semigroups arising from the proof of item (A2).
We
remark that an extension of item (A2) to the case of singular curves will
provide us with examples of numerical semigroups that cannot be realized as
Weierstrass semigroups even for singular curves.
\medskip

We have employed the methods used in [T], where one of the key
tools is  Castelnuovo's genus bound
for curves in projective space [C], [ACGH, p.116], [R, Corollary 2.8]. Since
the
coverings considered here can be of degree bigger than two, we will also
use as a key tool the other famous genus bound of Castelnuovo which concerns
subfields of the field of rational functions of the curve [C1], [St].

  \medskip

{\bf Conventions.} Throughout this paper, the word curve will mean a
projective, irreducible, non--singular algebraic curve defined over an
algebraically closed field $K$ of characteristic $p \ge 0$. By a
numerical
semigroup $H$ we will mean a subsemigroup of $(\bbN ,+)$ whose
complement
in $\bbN$ is finite. For $i \in \bbN$, we denote by $m_i = m_i(H)$ the
$i^{th}$ element of the semigroup $H$ and by $G(H)$ the set $\bbN \setminus
H$. When  $H$ is the Weierstrass semigroup of some point $P$, we
just write $H(P)$, $G(P)$ and $m_i = m_i(P)$. Given a curve $X$, the
symbols $g^{r}_{d}$, $K(X)$ and ${\rm div}_\infty (f)$ will denote
respectively an $r$-linear system of degree $d$ on $X$, the field of
regular functions of $X$ and the polar divisor of $f \in K(X)$.

\section{Semigroups of type $(N,\gamma )$.}

Let $H$ be a numerical semigroup. The natural number $g=g(H):=\#
G(H)$ is called the genus of $H$. The elements of $G(H)$
are called the gaps of $H$ and those of $H$ are called the non--gaps of $H$.

Fix a positive non-gap $m \in H$. For $i=1,\dots , m-1$, denote by
$s_i=s_i(H,m)$ the
smallest element of $H$ such that $s_i\equiv i$ (mod $m$) and then define
$e_i=e_i(H,m)$ by the equation
\begin{equation}
s_i=e_i\ m +i\ .
\end{equation}
By the semigroup property of $H$, we have that $e_i$ is the number of gaps
$\ell$ for which $\ell\equiv i$ (mod $m$). Consequently
\begin{equation}
g=\sum^{m-1}_{i=1}\ e_i
\end{equation}
and also
\begin{equation}
\left\{\matrix{ e_i+e_j\ge e_{i+j},\hfill &{\rm if}\ \  i+j<m.\cr
&\cr
e_i+e_j\ge e_{i+j-m}-1,&{\rm if}\ \ i+j>m.\cr}\right.     
\end{equation}
Conversely, given numbers $m,e_1,\dots e_{m-1}$ satisfying the above
relations one indeed has a semigroup. In particular, $m=m_1$ and the
respective $e_i$'s completely determine $H$ (cf. [H]). Let $N$ be a
positive integer. We associate to $H$ the number:
\begin{equation}
\gamma_N := \{\ell \in G(H): \ell \equiv 0\ ({\rm mod}\ N)\}.
\end{equation}

{\bf Lemma 2.1.} Let $H$ be a numerical semigroup of genus $g$, $N
\ge 2$ an integer. If $h\in H$ such that $gcd(h,N) = 1$, then
$$
h\ge {{2g-2N\gamma_N}\over {N-1}} +1\ .
$$

{\bf Proof.} Set $\gamma = \gamma_N$ and let $m=Nn$ be the least positive
non-gap of $H$ which
is  multiple of $N$. Then, $\gamma = \sum_{i=1}^{n-1}
e_{Ni}$ and  there exists $i\in  \{1,\ldots,m-1\}$ so that
$gcd(i,N)=1$ and $h\ge s_i$.

{\bf Claim.} For $k=1,\dots ,N-1;\ell =0,\dots ,n-1;ki+N\ell\not\equiv 0$
(mod $N$). Moreover these numbers are pairwise different modulo $Nn$.

Indeed, if $ki+N\ell\equiv 0$ (mod $N$), then $ki\equiv 0$ (mod $N$). Since
$gcd(i,N)=1$ we have $k\equiv 0$ (mod $N$), a contradiction. On the other hand,
$ki+N\ell\equiv k_1i+N\ell_1$ (mod $Nn$) implies $(k-k_1)i\equiv 0$ (mod $N$)
and so $k-k_1\equiv 0$ (mod $N$) which gives $k=k_1$. Consequently from
$N\ell\equiv N\ell_1$ (mod $Nn$) we obtain that $\ell\equiv\ell_1$ (mod $n$)
and so $\ell =\ell_1$.

For $k$ and $\ell$ as in the above claim, write $ki+N\ell =
a_{k\ell}Nn+r_{k\ell}$ (*) with $0<r_{k\ell}<Nn$. From (9) and induction on
$k$ and $\ell$ we have
$$
ke_i + e_{N\ell} \ge e_{r_{k\ell}} - a_{k\ell}
$$
\noindent where we assume $e_0:= 0$. Adding up these inequalities, from
the Claim and (8) we get
$$
{(N-1)Nn\ e_i\over 2} +(N-1)\gamma\ge g-\gamma -\sum_{k,\ell}\ a_{k\ell}\ .
$$
\noindent Now from (*) and the Claim we have
$\sum_{k,\ell}\ a_{k\ell}=(N-1)(i-1)/2$ and hence
the proof
follows from the above inequality and (7).\qquad $\Box$
\medskip

{\bf Remarks 2.2.} (i) The above lemma subsumes the following
result
due to Jenkins [J]: ``let $H$ be a numerical semigroup of genus $g$ and
$0<m<n$ non-gaps of $H$ so that $gcd(m,n) =1$; then $g \le
(m-1)(n-1)/2$ ".
Indeed, by using the notation of the lemma, take $N=n$; then $\gamma_N =0$
and Jenkins' result follows with $h =m$.\\
(ii) The lower bound of Lemma 2.1 is the number of ramified points minus
one of an $N$--sheeted covering of curves of genus $g$ and $\gamma_N$
respectively (defined over a field of characteristic $p\not\!\vert N$),
where all
the ramified points  are totally ramified points.\qquad $\Box$ \medskip

The next lemma will help us to understand the structure of the semigroups of
type $(N,\gamma)$.
\medskip

{\bf Lemma 2.3.} Let $H$ be a numerical semigroup.
\begin{list}
\setlenght{\rightmargin 0cm}{\leftmargin 0cm}
\itemsep=0.5pt
\item[(i)] Suppose that $H$
 fulfils
conditions (a) and (c) of Definition 2.
Set $F:=\{(2\gamma +i)N:i\in\bbZ^{+}\}$. Then:
\subitem(i.1) $F\subseteq H$, $2N\gamma \in H$,
\subitem(i.2) $\gamma = \gamma_N$.
\item[(ii)] Conversely, let $N>0$ be an integer. Then, $H$
fulfils condition (a) and (c) of  Definition 2 with $N$ and $\gamma_N$, and
$2N\gamma_N \in H$.
\end{list}

{\bf Proof.} (i) If $\gamma =0$ then $N\in H$ and so we have (i.1) and
(i.2). Let $\gamma\ge 1$ and denote by $f_1<\ldots <f_\gamma$ the
$\gamma$ positive multiples of $N$ non-gaps of $H$ in $[N,2N\gamma]$. So
$f_1 > N$. Suppose that $F\not\subseteq H$ and let $(2\gamma + i)N$ be the
least element of $F\cap
G(H)$. By the semigroup property of $H$ we have $(2\gamma +i)N-f_j\in
G(H)$ for
$j=1,\dots ,\gamma$. Then, by the selection of $(2\gamma +i)N$
we have that
$\{ (2\gamma+i)N-f_j:j=1,\dots ,\gamma\}$ are all the gaps $\ell$ of
$H$ such that
$\ell \equiv 0$ (mod $N$) and $\ell \le 2N\gamma$. The least of these
gaps satisfies $(2\gamma+i)N-f_\gamma \ge iN \ge 2N$ due to
condition (c) of Definition 2. Consequently $f_1 =N$ which is a
contradiction.
The statement (i.2) follows from (i.1) since the gaps $\ell$ for which
$\ell\equiv
0$ (mod $N$) belong to the interval $[N,2N\gamma ]$. It remains to
proof that $2N\gamma \in H$. Suppose that $f_\gamma <2N\gamma$. Then
we get $\gamma+1$ gaps multiples of $N$ namely,
$2N\gamma-f_\gamma,\ldots, 2N\gamma-f_1, 2N\gamma$ which is a
contradiction due to (i.2).\\
(ii) By the definition of $\gamma=\gamma_N$ (see (10)), there exist at least
$\gamma$
positive non-gaps - all of them being multiples of $N$ - in the interval
$[N,2N\gamma]$. Denote by $f_1< \ldots < f_\gamma$ such non-gaps. Let
$\ell$ be the biggest gap of $H$ so that $\ell \equiv 0$ (mod $N$). We claim
that $\ell < f_\gamma$, because on the contrary case we would have - as in
the previous
proof with $\ell$ instead of $2N\gamma$ - ($\gamma+1$) gaps which is a
contradiction with the definition of $\gamma$. This implies that
$2N\gamma, (2\gamma+1)N, \ldots $ are non-gaps and we are done.
\qquad $\Box$
\medskip

{\bf Corollary 2.4.} Let $H$ be a numerical semigroup, $\gamma$ a
non--negative integer, $M,N,r$ positive integers so that $2(\gamma
+r)M>(2\gamma +r)N$. Then $H$ cannot be both of type $(N,\gamma )$ and of type
$(M,\gamma +r)$.

{\bf Proof.} Suppose $H$ is both of type $(N,\gamma )$ and of type
$(M,\gamma +r)$. From the previous lemma and since $H$ is of type $(M,\gamma
+r)$ we have
$$
2(\gamma +r)M=m_{\gamma +r}\le (2\gamma +r)N\ .\qquad{\Box}
$$

Using Lemma 2.3 (ii) and Lemma 2.1 we have the following criteria for
the  type $(N,\gamma_N)$ of numerical semigroups.
\medskip

{\bf Corollary 2.5.} Let $H$ be a numerical semigroup and $N>0$ an
integer. Suppose that every $h \in H$ such that $h\not\equiv 0$ (mod $N$)
satisfies $h \ge 2N\gamma_N+1$. Then, $H$ is of
type  $(N,\gamma_N)$.\qquad $\Box$
\medskip

{\bf Corollary 2.6.} Let $H$ be a numerical semigroup of genus $g$ and
$N$  prime.
If $g>\rho_1(2\gamma_N,N,\gamma) =
N^2\gamma_N - N+1$, then $H$ is of type
$(N,\gamma_N)$.\qquad $\Box$
\medskip

We also have:
\medskip

{\bf Corollary 2.7.} Let $H$ be a semigroup of type $(N,\gamma)$
with $N$ prime. Let $A \ge \gamma +1$ be an integer and $g$ the genus of
$H$. If $g>\rho_1(A,N,\gamma)$ (see (1)), then
$$
gcd(m_1,\ldots, m_{A-\gamma}) = N.\qquad \Box
$$

\section{Proof of Theorem A.}

We study certain $N$--sheeted coverings
$$
\pi: X \to \tilde X
$$
\noindent of curves. To fix notation, we let $X$ and $\tilde X$ be
curves of genus
$g$ and $\gamma$, respectively. We assume that there is a point $P\in X$
such that $\pi$ is totally ramified at $P$, i.e., $X$ will be a curve of
type $(N,\gamma)$. We are mainly interested in relating the Weierstrass
semigroups at $P$ and $\tilde P:= \pi(P)$. Since $P$ is totally ramified,
 $\tilde
m_iN \in H(P)$ for $\tilde m_i \in H(\tilde P)$. Moreover, since $\tilde
m_{\gamma
+ j} = 2\gamma + j$  for $j \in \bbN$, we have the following statements:
\begin{list}
\setlength{\rightmargin 0cm}{\leftmargin 0cm}
\item[(I)] $\gamma_N = \gamma_N(P) := \#\{\ell \in G(P): \ell \equiv 0\
({\rm mod}\ N)\} \le \gamma.  $
\item[(II)] $ m_\gamma= m_\gamma(P) \le 2\gamma N \in H(P).$
\end{list}
\noindent
Note that equality in (II) implies equality in
(I), and,
 $H(P)$ is of type $(N,\gamma)$ if and only if equality in (II)
holds. Moreover, if $h\in H(P)$ so that $gcd(h,N)=1$, from Lemma 2.1
and (I) we  have
\begin{equation}
h \ge {{2g-2N\gamma_N}\over N-1} +1 \ge {{2g-2N\gamma}\over N-1}+1.
\end{equation}
\noindent Hence we have the following generalization of [T, Lemma 3.1].
\medskip

{\bf Lemma 3.1.} Assume the above notation and suppose $g >
\rho_1(2\gamma,N,\gamma)= N^2\gamma-N+1$. Then, every $h \in
H(P)$ such that  $$
h \le {g+N(N-2)\gamma \over N-1}
$$
\noindent satisfies $gcd(h,N) >1$. \qquad $\Box$
\medskip

The following result - due to Castelnuovo - will be used, among other
things, to prove the implication (ii) $\Rightarrow$ (i) of Theorem
A regardless of the characteristic of the base field and to construct
examples in
order to show that in some cases the  bounds of our results are sharp.
\medskip

{\bf Lemma 3.2 ([C1, St]).} Let $X$ be a curve of genus $g$ and
$K_{1}, K_{2}$ be subfields of $K(X)$ with compositum $K(X)$. If $n_{i}$
is the degree of $K(X)$ over $K_{i}$ and $g_{i}$ is the genus of $K_{i}$
for $i=1,2$, then
$$
g \le (n_{1}-1)(n_{2}-1) + n_{1}g_{1} + n_{2}g_{2}.\qquad \Box
$$

For  $N$ prime, we have the uniqueness of $\pi$ above
provided
$g$ is large enough, and we also have a criterion to decide when a point is
totally ramified:
\medskip

{\bf Corollary 3.3.} Let $X$ be a curve of genus $g$, $N$ prime and
$\gamma$ a non--negative integer.\\
(i) If
$$
g > \rho_5(N,\gamma):= 2N\gamma +(N-1)^2,
$$
\noindent then $X$
admits at most one $N$--sheeted covering of a curve of genus $\gamma$.\\
(ii) Let $P \in X$, $\tilde X$ be a curve of genus $\gamma$ and, $\pi$ an
$N$--sheeted covering map from $X$ to $\tilde X$. Then, $P$ is totally
ramified for $\pi$ provided there exists $h
\in H(P)$ such that
\begin{equation}
(N-1)h < g -N\gamma +N-1.
\end{equation}

{\bf Proof.} (i) If K(X) have two differents
fields $K_1$ and $K_2$ both of genus $\gamma$, then by Lemma 3.2 we have
$g \le \rho_5(N,\gamma)$.\\
(ii) Let $f\in K(X)$ with div$_\infty
(f)=hP$ and $K^{'}$ be the compositum of $K(f)$ and $K(\tilde X)$.
Using $N$ prime and the hypothesis on $h$, from Lemma 3.2 it follows that
$K^{'}=K(\tilde X)$. Then, there exists $\tilde f \in K(\tilde X)$ so
that $f = \tilde f\circ \pi$. Consequently,
  the ramification number of $\pi$ at $P$
is $N$ and so $P$ is totally ramified for $\pi$.
\rightline{$\Box$}

Next we look for conditions to have equality in (I) or (II).
\medskip

{\bf Lemma 3.4.} Let $X$, $\tilde X$, $\pi$, $P$, $\tilde
P$ and $N$ be as above. If either $p = char(K)$ $\not\!\vert N$ or  $N$
prime and $g >\rho_2(N,\gamma)$ (see (2)), then
$$
\gamma_N = \gamma.
$$

{\bf Proof.} It will be enough to show
that: $nN \in H(P)
\Rightarrow n \in H(\tilde P)$.\\
{\bf Case 1:} $p \not\!\vert N$. Let $z$ be a local parameter at $P$ so
that $z^N$ is also a local parameter at $\tilde P$. Let $\Psi$ (resp.
$\tilde \Psi $)
denote the inmersion of $K(X)$ (resp. $K(\tilde X)$) into the field of
Puiseux series at
$P$ (resp. $ \tilde P$) $F_1 = K((z))$ (resp. $F_2 = K((z^N))$). Since
${\Psi \mid}
_{K(\tilde X)} =  \tilde \Psi$ we have that $Tr_{F_1\mid F_2}\circ
\Psi = \tilde \Psi \circ Tr_{K(X)\mid K(\tilde X)}$ (*) (Tr means trace).
Let $f \in K(X)$ with ${\rm div}_\infty (f) = nN$. Write $f =
\sum_{i=-nN}^{\infty} c_iz^i$. Then, by considering the base
$\{1,z,\ldots , z^{N-1} \}$ of $F_1\mid F_2$, we have that $Tr_{F_1\mid
F_2} (f) = \sum_{i=-n}^{\infty} {Nc_{iN}z^{iN}}$. Consequently, from
(*) it follows that
the order of $\tilde f:= Tr_{K(X)\mid K(\tilde X)}(f)$ at $\tilde P$ is $n$
and, since
$f$ has no other pole, ${\rm div}_\infty(\tilde
f) = n\tilde P$ and we are done.\\
{\bf Case 2:} $g > \rho_2(N,\gamma)$. From the proof of Corollary 3.3
(ii) we have that
$f = \tilde f \circ \pi$ for some $\tilde f \in K(\tilde X)$ whenever
${\rm div}_\infty(f) = hP$ with $h$ satisfying (12). Now, from the
hypothesis on $g$ we can applied the above statement for $h \in H(P)$ with
$h \le 2N\gamma - N$.  \qquad $\Box$
\medskip

{}From (I), (II), (11) and the lemma above we obtain:
\medskip

{\bf Corollary 3.5.} Assume the hypothesis of Lemma 3.4.\\
(i) H(P) satisfies conditions (a), (c) of Definition 2 (with $N$ and
$\gamma$) and $2N\gamma \in H(P)$.\\
(ii) Suppose $N$ is prime. If either $N \not= p$ and $g >
\rho_1(2\gamma,N,\gamma) = N^2\gamma -N+1$ or
 $g > \rho_2(N,\gamma) $, then $H(P)$ is of type $(N,\gamma)$.\qquad
$\Box$
\medskip

{\bf Remark 3.6.} Let $\pi: X \to \tilde X$ be an $N$--sheeted
covering of curves of genus $g$ and $\gamma$ respectively. Assume $g >
\rho_2(N,\gamma)$ and hence, in particular that $\pi $ is ``strongly
branched" (cf.
[A]). When $\pi$ is a ``maximal strongly branched" (e.g. $N$ prime) we
still have the result in Lemma 3.4 [A, Lemma 4].\qquad $\Box$.
\medskip

To deal with the ``geometry" of Theorem A we need the
other Castelnuovo genus bound lemma:
\medskip

{\bf Lemma 3.7 ([C], [ACGH, p.116], [Ra, Corollary 2.8]).} Let $X$ be a
curve of genus $g$ that
admits a birational morphism onto a non--degenerate curve of degree $d$ in
$\bbP^r(K)$. Then
$$
g\le c(d,r):={m(m-1)\over 2}\ (r-1)+m\varepsilon
$$
where $m$ is the biggest integer $\le (d-1)/(r-1)$ and
$\varepsilon = d-1 - m(r-1)$. \rightline{$\Box$}

{\bf Lemma 3.8.} Let $X$ be a curve of genus $g$, $N$ a prime and
$\gamma \ge 0$ an integer. Let $A\ge 2\gamma +1$ be an integer satisfying
the hypotheses (5) and (6) of  item (A4) ( Theorem A).
If $X$ admits a base--point--free linear system  $g^{A-\gamma}_{AN}$,
then $X$ is an $N$--sheeted covering of a curve of genus $\gamma$.

{\bf Proof.} Let $\pi :X\to \bbP^{A-\gamma}(K)$ be the morphism defined
by $g^{A-\gamma}_{AN}$.
\medskip

{\bf Claim:}\quad $\pi$ is not birational.

If by way of contradiction $\pi$ is birational, we can applied the lemma
above to obtain $g \le c(AN,A\gamma)=\rho_4(A,u,N,\gamma)$.
\smallskip

Let $t$ be the degree of $\pi$ and
$\tilde X$ the normalization of $\pi (X)$. Then the induced
morphism
$\pi :X\to\tilde X$ is a covering map of degree $t$ and $\tilde
X$ admits a base--point--free linear system $\tilde g^{A-\gamma}_
{AN\over t}$. In particular we have $t \le AN/(A-\gamma)$ and the
hypothesis (5) implies $t=N$. Now,
by the Clifford's theorem we have that $\tilde g^{A-\gamma}_{A}$ is
nonspecial, and consequently by the Riemann--Roch theorem the genus
of $\widetilde X$ is $\gamma$ and the proof is complete.\qquad $\Box$
\medskip

{\bf Proof of Theorem A.} {\bf (A1):} Corollary
3.5.\\
{\bf (A2):} Since $\rho_3(N,\gamma)> \rho_1(2\gamma+2,N,\gamma)$ from
Corollary 2.7 we have $D:= gcd(m_1(P), \ldots, m_{\gamma+2}) = N$. In
particular $m_{\gamma + 2} = (2\gamma +2)N$.
Now, we can apply the Claim in the proof above with $A=2\gamma+2$ and
$\rho_4(2\gamma+2,N-1,N,\gamma)= \rho_3(N,\gamma)$ to conclude that
the degree $t$ of the rational map obtained from the liner system
$|m_{\gamma
+2}P|$ is bigger than 1. Due to the fact that $t\vert D$ and $N$ prime,
we  conclude  that $t=N$, and by a similar argument to the above
proof (last lines)  we see that the covered curve has genus $\gamma$ and we
are done.\\
{\bf (A3):} By Corollary 2.7 we have $m_{A-\gamma}(P) = AN$ and it
follows the proof.\\
{\bf (A4):} The  above lemma shows that $X$ is an $N$--covering of a curve
of genus $\gamma$. Since the covering is given by $|ANP|$, we have that
$P$ is a totally ramified point of $\pi$ and the proof is complete. \qquad
$\Box$ \medskip

{\bf Corollary 3.9.} (i) Let $\pi :X\to\tilde X$
be an $N$--sheeted covering of curves of genus $g$ and $\gamma$ respectively.
Suppose $N$ prime and $g> \rho_3(N,\gamma)$. Then $P$ is
totally ramified for
$\pi$ if and only if $H(P)$ is a semigroup of type $(N,\gamma )$.\\
(ii) Let $H$ be a Weierstrass semigroup of genus $g> \rho_3(N,\gamma)$.
 Then $H$ is of type $(N,\gamma )$ if and only if there exists
an integer $A \in [2\gamma +2, 2\gamma+2 +{\gamma \over N-1}]$
satisfying (5) and such that $m_{A-\gamma} =AN$.

{\bf Proof.} (i) Proof of item (A2) of Theorem A.\\
(ii) For the numbers $A$ in that interval we have
$\rho_1(A,N,\gamma) \le {N(2N+1)\gamma \over 2} + (N-1)^2 <
\rho_3(N,\gamma)$ and hence the ``if" part of the statement is just item
(A3). To prove the ``only if" part notice that $u(A) = N-1$ and hence
$\rho_4(A,u(A),N,\gamma) = \rho_3(N,\gamma)$. Now the hypotheses on $A$
assure
that the degree of the map obtained from $|ANP|$ is $N$, and since this
number is the $g.c.d$ of the non--gaps $m_1,\ldots,m_{A-\gamma}$, it follows
the proof. \qquad $\Box$
\medskip

{\bf Remark 3.10.} Remark 3.11 (ii), (iii) below show that neither the
bound on $g$ nor the hypothesis (5) of the above corollary (part (ii)) can be
dropped. It
would be interesting to have an arithmetical proof  of Corollary 3.9 (ii)
(i.e. without the assumption that $H$ is a Weierstrass semigroup),
because any  counter example  to the above
question would be a numerical semigroup that cannot be realized as
Weiertrass semigroup. This numerical semigroup could be see as a ``mid
term" between the examples stated in the last section. The most simple case of
the above question is
for $N=2$. But this case does not provide any counter example [T1].
\medskip

{\bf Remarks 3.11.} Let $N$, $\gamma$ be a prime and a non--negative
integer respectively and suppose $p \not\!\vert 2N$.\\
{\bf (i)} The following example has respect to
Lemma 2.1, Corollary 2.7, Lemma 3.1, Corollary 3.5 (ii) and item
(A1) of Theorem A.
Let $g>0$ be an integer such
that $g-N\gamma \equiv 0$ (mod $(N-1)$) and $L := {{2g-2\gamma N}\over
{N-1}}+1$ is coprime with $2N$. Define $i_1:= 2\gamma+1$ and
$i_2:= {g-(2N-1)\gamma \over N-1} $ (hence $i_1 +2i_2 =L$). For $j=1,\ldots,
i_1, k=1, \ldots i_2$, choose $a_j, b_k$  pairwise
distinct elements of $K$. Now consider the curve $X$ defined by the
equation
\vspace{-0.5pt}
$$ y^{2N} = \mathop{\Pi}\limits^{i_1}_{j=1}\ (x-a_j)
\mathop{\Pi}\limits^{i_2}_{k=1}\
(x-b_k)^2
$$
\vspace{-0.5pt}
Then, by the Riemann--Hurwitz relation  we have that
 the genus
of $X$ is $g$. Moreover, $X$ is a $N$--sheeted covering of the
hyperelliptic curve $\tilde X$ of genus $\gamma$ whose field of rational
functions is $K(x,z)$,
 where $z = y^{N}/\Pi_{k=1}^{i_2}\ (x-b_k)$. Since
$gcd(L,2N)=1$, there exist just one point $P\in X$ over $x = \infty$ and
consequently $X$ is a curve of type $(N,\gamma)$ over $\tilde X$.
\medskip

\noindent {\bf Claim :}\qquad $H(P) = H:= \langle 2N, L, (2\gamma +1)N
\rangle$. \medskip

\noindent This claim shows that the result in Lemma 2.1 is the best
possible. Considering $g = \rho_4(A,N,\gamma)$, we have $L = AN-1$ and
hence $m_{A-\gamma}(H) \le AN-1 $ provided $A \ge 2\gamma$. Hence the claim
also shows the sharpness of the bound on $g$ of Corollary 2.7. By
specializing $A = 2\gamma$ we also see the sharpness of the bound on $g$
of Lemma 3.1 and Corollary 3.5 (ii) (case $p \not\!\vert N$) respectively.
With respect to item (A1) is not difficult to see that in the above curve
$X$
(with $A = 2\gamma$), all the Weierstrass semigroups at totally ramified
points are not of type $(N,\gamma)$ (the case $N=2$ is in [T, Remark
3.9]). However, we cannot say the same about the other Weierstrass
semigroups of $X$.
\medskip

\noindent {\bf Proof of the Claim.} Since the genus of $H$ is at most
$g$, it will
be enough to show that $H \subseteq H(P)$. This is true because, ${\rm
div}_\infty (x) = 2N P$, ${\rm div}_\infty (y) = L P$, and $(2\gamma +1)N
\in H(P)$ due to the fact that $P$ is totally ramified over $\tilde X$
which has genus  $\gamma$.\qquad $\Box$
\medskip

\noindent {\bf (ii)} This example is related to the bound on the genus in
Lemma  3.8 and item (A4) of Theorem
A. Set $i_1 := 2N\gamma +2N-1$. The curve $X$ defined by the equation
\vspace{-0.5pt}
$$
y^{2N} = \mathop{\Pi}\limits^{i_1}_{j=1}\ (x-a_j),
$$
\vspace{-0.5pt}
\noindent where the $a_j's$ are pairwise distinct elements of $K$, has the
following properties:
\begin{list}
\setlength{\rightmargin 0cm}{\leftmargin 0cm}
\itemsep=0.5pt
\item[(1)] its genus is $g = N(2N-1)\gamma + (N-1)(2N-1) = \rho_3(N,\gamma)$;
\item[(2)] the Weierstrass semigroup at the unique point $P$ over
$x = \infty$ is generated by $2N$ and $i_1$;
\item[(3)] $m_{A-\gamma}(P) = AN$ provided $2\gamma \le A
< 4\gamma +4 - {2\over N}$;
\item[(4)] it cannot be an $N$--covering of a curve of genus $\gamma$.
\end{list}
\noindent Consequently, the upper bound $\rho_3(N,\gamma)$ for the genus
in both Lemma 3.8 and item (A4) is necessary for $2\gamma +1 \le A \le
4\gamma +4 -{2\over N}$. In particular, if N-1 is the
biggest integer $\le {N\gamma +N-1 \over A-\gamma -1}$, $\rho_3(N,\gamma)$
is sharp. In the other cases, we don't know the sharpness of
$\rho_4(A,u,N,\gamma)$.
\medskip

\noindent {\bf Proof of properties (1)--(4).} (1) follows from
Riemann--Hurwitz
relation. To prove (2), we notice that ${\rm div}(x) = 2NP$ and
${\rm div}(y) = i_1P$ and so $H(P) \supseteq \langle 2N,
i_1 \rangle $. Since the last semigroup also has genus $g$, we have (2).
To prove (3) we notice that in the interval $[1,AN]$ the number of
multiples of $N$ non-gaps
 of $H(P)$ is $A\over2$ (or $(A+1)/2$);
it
has ${A\over2} -\gamma$ (or $(A-1)/2 -\gamma $) non--gaps which are
congruent to $2N-1$ module
$N$, and its other non--gaps are bigger than $AN$  (here we use $A <
4\gamma +4 - {2\over N}$). Finally, if $X$ is an $N$--sheeted covering of a
curve of
genus $\gamma$, by Castelnuovo's lemma (Lemma 3.2) the genus $g$ would be
at most $(2N\gamma +2N-2)(N-1) + N\gamma < \rho_3(N,\gamma)$, a
contradiction.\\
{\bf (iii)} Here, we show that the arithmetical conditions
(5) cannot be dropped if we suppose
\vspace{-0.5pt}
$$
g > {\rm max}\ \{{1\over 2}AN[A(N-2) +2\gamma +3],
A(N-1)(N-2) +(3N-2)\gamma +3(N-1)\}.
$$

Since $A\ge 2\gamma +1$ we have ${AN\over A-\gamma} \le 2N-1$
and then one has to check (5) among the integers of
the set $[2,2N-1] \setminus {N}$. Let $t$ be an integer of the above set
such that $t\vert A$ and set $r:= {AN\over t} -A +\gamma +1$. Let $g>0$
be an integer satisfying the above bound and such that $i_1:= {2g \over
rt -1} +1 $ is also an integer. The curve in the previous remark, with
$rt$ instead of $2N$ and the above $i_1$, has genus $g$ and just one point
$P$
over $x=\infty$ which satisfies $m_{A-\gamma}(P) = AN$ (here we use
the first
part of the bound). But $X$ cannot be an $N$--sheeted covering because
on the contrary by Castelnuovo's lemma (Lemma 3.2) the genus would be at
most the second part of the above bound.\\
{\bf (iv)} Finally, some words about items (A2) and (A3). Since statement (ii)
of  Theorem A is stronger than (iii), one might expect to sharpen
$\rho_3(N,\gamma)$ (this would be relevant to the examples in the next
section). In order to do that, one might use
 Castelnuovo's
theory (cf. e.g. [E-H-1, $\oint$ 3], [Ci]) or ``results" extending this
theory to Hilbert functions of points in projective spaces [E-G-H].
Specifically, one could use analogous bounds to $c(d,r)$ in order to deal
with  curves of  genus $g \le
\rho_3(N,\gamma) = c(AN,A-\gamma)$. The point is that
one knows how must look the curves whose genus attain the
mentioned bounds. For instance,
one  can applied the above considerations to double covering of curves of
genus one or two and the result is that (A2) is still valid for
$g \ge \rho_3(2,\gamma)-2$ ($\gamma \in \{1,2\}$) (see also [G, Lemmas
7 and 9]). In general, we think that item (A2) must be true with a bound of
type ``$\rho_3(N,\gamma) - N$". We remark that by applying the arithmetical
properties of
 semigroups of type $(N,\gamma)$ one can find a ``kind of algorithm to
compute Hilbert functions". We will intend to describe this in a later
paper.

With respect to the sharpness of the bound on $g$ of item (A3), we just
want to say that it depends on the existence of certain Weierstrass
semigroups.\qquad $\Box$
\section{Hurwitz's question.}
In this section we construct numerical semigroups with $\ell_g$ given,
that cannot be realized as Weierstrass semigroup.
These examples will include symmetric and
quasi--symmetric semigroups generalizing those in [T, Scholium 3.5] and
[O-S, Example 6.5].
\vspace{-0.5pt}
\subsection{Corollaries of Buchweitz's criterion; case $n=2$.}
Let $X$ be a curve, $H$ a numerical semigroup both of genus $g$.
Denote by $\ell_1=\ell_1(H) <,\ldots,\ell_g= \ell_g(H)$ (resp.
$G_n=G(H)$) the gaps (resp. the set of
sums of $n$ gaps) of $H$. By the definition of $g$ and by the semigroup
property of $H$, we have $g\le \ell_g \le 2g-1$ ([B, O]). Semigroups with
$\ell_g = g$ are realized for all but finitely many points of $X$
provided it is defined in characteristic zero or characteristic larger
than 2g-2 (see e.g. [S-V]). In the remaining cases the situation can be
different (see e.g. [Sch], [G-V]). On the other hand, if $\ell_g = 2g-1$,
$H$ is called symmetric because between the non--gaps and gaps of $H$ we
have the following property: $S(1)$: ``$h \in H \Leftrightarrow \ell_g -h
\in
G(H)$". These semigroups are important at least for two reasons : 1)
they arise in a natural way in the context of Gorenstein rings (cf.
[Ku]); 2)
for $n \ge 2$, $\#G_n = (2n-1)(g-1)$ ([O, Thm. 1.5]), i.e., they satisfy
Buchweitz's criterion (BC) (see Section 1).

What  can be said if $\ell_g < 2g-1$ ?. If $\ell_g = 2g-2$, $H$ still
satisfies  property $S(1)$ except for $g-1$ [O, Prop. 1.2] ( consequently
$H$ is called quasi--symmetric).  Then, it is not surprising that
quasi--symmetric semigroups also
satisfy (BC). In fact, here we have $\# G_n = (2n-1)(g-1) -(n-2)$
([O-S, Thm. 1.1]). As in the case of the symmetric semigroups, these
semigroups can also be realized as Weierstrass semigroups of
Gorestein curve, but in view of the term ``$-(n-2)$" above, here one has to
allow reducible curves
 (cf. [O-S, $\oint$ 3]). In general we have the following properties of
type $S(1)$. Suppose $\ell_g \in \{2g-2i+1, 2g-2i\}$ with $i\ge 1$.
Considering the pairs $(r, \ell_g -r)$ for $r =1,\ldots, g-i $, $H$
satisfies the property
\medskip

\noindent $S(i):$\quad If $\ell_g$ is odd, $H$ fulfils property $S(1)$
except $2i-2$ gaps of type: $g-i < h_{i-1}<\ldots < h_1 <\ell_g$, $\ell_g -
h_{i-1}, \ldots, \ell_g - h_1$. If $\ell_g$ is even, $H$ fulfils
condition $S(1)$ except $2i-2$ gaps of the above type and the gap $(g-i)$.
\medskip

 Since for $i>1$ we have gaps different from those arising
in $S(1)$, it
seems to be difficult to obtain a  closed form for $\#G_n$. However, we
think that the following must be true provided $\ell_g \le 2g-2$:
\vspace{-0.5pt}
$$
G_n = \{n,n+1,\ldots,(n-1)\ell_g \}
\mathop{\cup}\limits^{g}_{k=1}\{(n-1)\ell_k +\ell_j : j=1, \ldots g \}.
$$
\vspace{-0.5pt}
{}From the proof of [O-S, Thm. 1.1] at least the inclusion
``$\supseteq$'' holds. In particular we have
\vspace{-0.5pt}
$$
\# G_2 = (\ell_g -1) + g + \Lambda,
$$
\vspace{-0.5pt}
 where $\Lambda$ is a non--negative integer. Consequently if
$\Lambda \ge 2i-2$ (resp. $2i-1$) for $\ell_g = 2g-2i+1$ (resp. $2g-2i$),
by Buchweitz's criterion $H$ is not a Weierstrass semigroup. Now,
consider the following sequences of gaps obtained from
those of property $S(i)$:
\vspace{-0.5pt}
$$
2h_1 >\ldots > h_1 +h_{i-1},\ \ {\rm and}
$$
\vspace{-0.5pt}
$$
2h_2>\ldots >2h_{i-1}.
$$
\vspace{-0.5pt}
 With the above notation we have
\medskip

{\bf Lemma 4.1.1.} Let $H$ be a numerical semigroup with $\ell_g =
2g-2i+1$ and $i \ge 4$. If $h_1 + h_{i-1} > 2h_2$ and $2\ell_g - h_u - h_v
\in G(H)$ for $(u,v)= (1,1),\ldots, (1,i-1)$ and $(u,v) =
(2,2),\ldots,(2,i-1)$, then $H$ is not a Weierstrass semigroup.

{\bf Proof.} The hypothesis involving $G(H)$ means that $h_u + h_v$ is not
the sum of $\ell_g$ with some other gap. Consequently from the sequences
of gaps above we have $\Lambda \ge 2i-2$ and it follows the proof.\qquad
$\Box$
\medskip

Using this criterion we can exhibit numerical semigroups  with a fixed
last gap which cannot be realized as Weierstrass semigroups. The
following example with $i=4$, $g=16$ is the well known Buchweitz's
semigroup.
\medskip

{\bf Corollary 4.1.2.} Let $g, i$ be integers so that $g \ge 9i -20$,
$i \ge 4$ and $3g + 5i -20 $ even, say equal to $2h_1$. Then the
numerical semigroup whose gaps are
$$
\{1,2,\ldots,g-i, h_1 - (a+2(i-3)),\ldots, h_1-(a+2), h_1 -a, h_1,
2g-2i+1 \},
$$
\noindent where $a= 2i-5$ is not a Weierstrass semigroup.\qquad $\Box$
\medskip

In the above examples one can use $a > 2i-6$ provided that $g \ge 2a -10
+5i$ and $3g+2a +i -10$ even.
\subsection{An application of item (A2) of Theorem A.}
First we notice that from the proof of item (A2), if $N$ is prime
and $H = \{m_0=0,m_1,\ldots
\}$ is a Weierstrass semigroup of type $(N,\gamma )$ of genus
$g>\rho_3(N,\gamma)$, then
the numerical semigroup
$$
\pi( H) = \{{m_i\over N}:1\le i\le\gamma\}\cup\{
2\gamma +i:i\in\bbN\}
$$
\noindent is also a Weierstrass semigroup.
We use this remark to prove an analogue of Corollary 4.1.2. The
semigroups of this result are also inspired by the properties $S(i)$ of
the last subsection.
Fix a numerical semigroup $\tilde H$ of genus $\gamma$ such that it is
not a Weierstrass semigroup.
Let $N$ be a prime and $g$ and integer. Write $g = \lambda N +u$ with
$0\le u <N$. Let $f$ be an integer such that $f \le u$ if $u>0$ and $f
<N$ otherwise. Set $N\tilde H:= \{hN; h \in \tilde H\}$. We are only going
to consider the case $2g-f \not\equiv 0$ (mod $N$) because in the other
case we can replace $g$ by $g+1$.
\medskip

{\bf Corollary 4.2.1.} With the above notation, consider the following sets
\begin{list}
\setlenght{\rightmargin 0cm}{\leftmargin 0cm}
\itemsep=0.5pt
\item[(1)]
$H_1= N\tilde H\cup\{ 2g-f-r:r\not\in N\tilde H,r\le
g-1\}$, if $2u \not\in [N,f+N]$
\item[(2)] $H_2 = H_1 \setminus \{e\}$ otherwise; where $e$ is the
biggest integer $\le (2g-f)/2 $.
\end{list}

If $g > \rho_3(N,\gamma)$, then $H_1$ and $H_2$ are
 numerical semigroups of type $(N,\gamma )$ of genus $g$
 whose last gap is $2g-f$
which are not
Weierstrass semigroups.
\medskip

{\bf Proof.} Set $H$ for $H_1$ or $H_2$. We notice that $\pi(H) =
\tilde H$ and so it will be enough to prove the arithmetical statements.
By the definition of $H$ and the hypothesis
on $g$, it follows that $H$ is a
semigroup of type $(N,\gamma )$ such that $H\supseteq\{ 2g,2g+1,\dots\}$ and
$\ell_g(H)=2g-f$ (here we use $2g-f\not\equiv 0$ (mod $N$)).
Consider $H=H_1$ and let
\begin{eqnarray*}
U & = & \#\{ h\in H:h\le 2g,h=2g\ \ {\rm or}\ \ h\equiv 0\ ({\rm mod}\
N)\},\\
V & = & \#\{ h\in H:h<2g,h\not\equiv 0\ ({\rm mod}\ N)\}\ .
\end{eqnarray*}
Then
\vspace{-0.5pt}
$$
U = \left\{\matrix{
2\lambda -\gamma,\hfill &{\rm if}\ \  u=0;\hfill\cr
2\lambda +1-\gamma ,\hfill&{\rm if}\ \  0<2u<N;\cr
2\lambda +2-\gamma ,\hfill &{\rm if}\ \  2u>N;\hfill\cr}\right.
$$
\vspace{-0.5pt}
$$
V=g-1-\#\{ 1\leq r\leq g-1:r\in N\widetilde H\}-\#\{ 1\leq r\leq g-1:r\equiv
2u-f ({\rm mod}\ N)\}\ .
$$
\noindent Since $2u-f\not= N$, we have
$$
V=\left\{\matrix{
g-1-(\lambda -1-\gamma )-\lambda ,\hfill &{\rm if}\ \  u=0;\hfill\cr
g-1-(\lambda -\gamma )-\lambda ,\hfill &{\rm if}\ \  2u-f<N;\cr
g-1-(\lambda -\gamma )-(\lambda +1),\hfill &{\rm if}\ \
2u-f>N.\hfill\cr}\right.
$$
Then $U+V=g$ is the number of non--gaps $\le 2g$ of $H_1$ because $2u
\not\in [N,n+f]$. In the other case, from the above computations we get $U+V
= g+1$ and since we have excluded $e$ we are done.\qquad $\Box$
\medskip

The last two corollaries give us numerical semigroups arising from
different phenomena. Moreover, notice that the ``intersection " of both
families of examples is empty. By considering the distribution of the
respective gaps sequences, we can think about of these examples as being
the ``extremal" cases of numerical semigroups that cannot be realized as
Weierstrass semigroups. Finally,
we would like to know if any
numerical semigroup which is not a Weierstrass semigroup but satisfies
Buchweitz's criterion must be of type $(N,\gamma)$ for some $N$ and $\gamma$.
\bigskip

\centerline{\bf ACKNOWLEDGMENTS}
\medskip

Thank you very much to the International Atomic
Energy Agency and UNESCO for the hospitality at the International Centre for
Theoretical Physics, Trieste. Special thanks to Profs. Oliveira
and St\"ohr who showed me how to applied item (A2) of Theorem A to the case
of quasi-symmetric
semigroups (in fact, their paper [O-S] suggested many ideas to this work);
Prof. Ulrich for pointing out
to me the remark on semigroups rings stated in the introduction and the
referee for  suggesting to improve the earlier version of the paper.
\bigskip

\centerline{\bf REFERENCES}
\small
\begin{description}
\itemsep=-0.5pt
\item{[A]}
Accola, R.D.M.: Strongly branched coverings of closed Riemann Surfaces, Proc.
Amer. Math. Soc. {\bf 26}, 315--322 (1970).

\item{[ACHG]}
Arbarello, E., Cornalba, M., Griffiths, P.A. and Harris, J.: Geometry of
algebraic curves, Vol.I, Springer--Verlag, New York (1985).

\item{[B]}
Buchweitz, R.O.: \"Uber deformationem monomialer kurvensingularit\"aten und
Weierstrasspunkte auf Riemannschen fl\"achen, Thesis, Hannover (1976).

\item{[B1]}
Buchweitz, R.O.: On Zariski's criterion for equisingularity and
non--smoothable monomial curves, Th\'ese, Paris VII (1981).

\item{[C]}
Castelnuovo, G.: Ricerche di geometria sulle curve algebriche, Atti. R. Acad.
Sci. Torino {\bf 24}, 196--223 (1889).

\item{[C1]}
Castelnuovo, G.: Sulle serie albegriche di gruppi di punti appartenenti ad una
curvea algebrica, Rendiconti della Reale Accademia dei Lincei (5) {\bf 15},
337--344 (1906).

\item{[Ci]}
Ciliberto, C.: Hilbert functions of finite sets of points and the genus
of a curve in a projective space, Space Curves, Rocca di Papa 1985, ed. F.
Ghione, C. Peskine and E. Sernesi (1987).

\item{[E-G-H]}
Eisenbud, D., Green, M., Harris, J.: Some Conjectures Extending
Castelnuovo Theory, preprint (1993).

\item{[E-H]}
Eisenbud, D., Harris, J.: Existence, decomposition and limits of certain
Weierstrass points, Invent. math. {\bf 87}, 499--515 (1987).

\item{[E-H-1]}
Eisenbud, D., Harris, J.: Curves in projective space, Les Preses de
l'Universit\'e de Montr\'eal, Montr\'eal, (1982).

\item{[G]}
Garcia, A.: Weights of Weierstrass points in double coverings of curves
of genus one or two, Manuscripta Math. {\bf 55}, 419--432 (1986).

\item{[G-V]}
Garcia, A., Viana, P.: Weierstras points on certain non--classical
curves, Arch. Math. {\bf 46}, 315--322 (1986).

\item{[H]}
Hurwitz, A.: \"Uber algebraische gebilde nit eindeutigen transformationen in
sich, Math. Ann. {\bf 41}, 403--441 (1893).

\item{[J]}
Jenkins, J.A.: Some remarks on Weierstrass points, Proc. Amer. Math. {\bf 44},
121--122 (1974).

\item{[K]}
Kato, T.: On the order of a zero of the theta function, Kodai Math. Sem. Rep.
{\bf 28}, 390--407 (1977).

\item{[Ku]}
Kunz, E.: The value-semigroup of a one-dimensional Gorenstein ring, Proc.
Amer. Math. Soc. {\bf 25}, 748--751 (1970).

\item{[O]}
Oliveira, G.: Weierstrass semigroups and the canonical ideal of non--trigonal
curves, Manuscripta Math. {\bf 71}, 431--450 (1991).

\item{[O-S]}
Oliveira, G., St\"ohr, K.O.: Gorenstein curves with quasi-symmetric
Weierstrass semigroups, Preprint (1994).

\item{[P]}
Pinkham, H.: Deformations of Algebraic Varieties with $\bbG_m$ Action,
Asterisque {\bf 20}, Soc. Math. France, (1974).

\item{[R]}
Rathmann, J.: The uniform position principle for curves in characteristic
$p$, Math. Ann. {\bf 276}, 565--579 (1987).

\item{[R\"o]}
R\"ohrl, H.: Unbounded coverings of Riemann surfaces and extensions of rings of
meromorphic functions, Trans. Amer. Math. Soc. {\bf 107}, 320--346 (1963).

\item{[Sch]}
Schmidt, F.K.: Zur arithmetischen Theorie der algebraischen Funktionen,
II, Allgemaire Theorie der Weierstrasspunke, Math. Z. {\bf 45}, 75--96
(1939).

\item{[S]}
St\"ohr, K.O.: On the modulli spaces of Gorenstein curves with symmetric
Weierstrass semigroups, J. reine angew. Math. {\bf 441}, 189--213 (1993).

\item{[S1]}
St\"ohr, K.O.: On the poles of regular differentials of singular curves,
Bol. Soc. Bras. Mat. {\bf 24}, 105--136 (1993).

\item{[S-V]}
St\"ohr, K.O., Voloch, J.F.: Weierstrass points and curves over finite
fields, Proc. London Math. Soc. (3), {\bf 52}, 1--19, (1986).

\item{[St]}
Stichtenoth, H.: Die ungleichung von Castelnuovo, J. Reine Angew. Math. {\bf
348}, 197--202 (1984).

\item{[T]}
Torres, F.: Weierstrass points and double coverings of curves with application:
Symmetric numerical semigroups which cannot be realized as Weierstrass
semigroups,  Manuscripta Math. {\bf 83}, 39--58 (1994).

\item{[T1]}
Torres, F.: Remarks on Weierstrass semigroups at totally ramified points
(provisory title), in preparation.
\end{description}

\end{document}